\documentclass[11pt, reqno]{amsart}
\usepackage{amssymb,latexsym,amsmath,amsfonts, amsthm}

\usepackage{latexsym}

\usepackage{color}
\usepackage[mathscr]{eucal}
\usepackage{comment}

\usepackage[linkcolor=blue, citecolor=black]{hyperref}
\hypersetup{colorlinks=true, urlcolor=blue}
\usepackage[nameinlink]{cleveref}
\usepackage{enumerate}

\numberwithin{equation}{section}

\theoremstyle{definition}
\newtheorem{definition}{Definition}[section]
\newtheorem{example}[definition]{Example}

\theoremstyle{remark}
\newtheorem{remark}[definition]{Remark}
\theoremstyle{plain}
\newtheorem{theorem}[definition]{Theorem}
\newtheorem{lemma}[definition]{Lemma}
\newtheorem{proposition}[definition]{Proposition}

\newtheorem{result}[definition]{Result}

\newtheorem{question}[definition]{Question}

\newcommand{\poly}{\mathscr{P}}

\newcommand{\C}{\mathbb C}
\newcommand{\R}{\mathbb R}

\begin{document}
	\title[Totally real discs]{Polynomial Convexity and Non-Polynomial Convexity of Totally Real Discs in $\mathbb C^2$}
	\author{Golam Mostafa Mondal}

	\address{Department of Mathematics, Indian Institute of Science, Bengaluru-- 560012, India}
	
	\email{golammostafaa@gmail.com, golammondal@iisc.ac.in}
	
	\keywords{Polynomial convexity; Totally real discs; Strongly pseudoconvex domains; Levi-flat hypersurfaces; Runge domains}
	\subjclass[2020]{Primary: 32E20, 32T15, 32V40}

\begin{abstract}
We study the polynomial convexity of smooth totally real
discs contained in the boundaries of strongly pseudoconvex
domains in $\mathbb{C}^2$. Using Wermer's classical example, we
construct a bounded domain with smooth strongly
pseudoconvex boundary containing a non-polynomially convex
totally real disc. We also construct a nonsingular
Levi-flat hypersurface in a non-Runge domain containing
such a disc. 
This provides a partial negative answer, to a conjecture proposed in an earlier work. Furthermore, we prove that the Wermer graph over the closed
disc of radius $R$ is polynomially convex for every
$R\in (0,1)$. Finally, we establish a sufficient condition
for the polynomial convexity of totally real discs contained
in the smoothly bounded strongly
pseudoconvex domains in $\mathbb{C}^2$.
\end{abstract}

\maketitle

\section{Introduction}
\label{sec:introduction}

Let $K$ be a compact subset of $\mathbb{C}^n$. The polynomial
convex hull of $K$ is defined by $\widehat{K}
=
\left\{
z\in\mathbb{C}^n:
|P(z)|\leq\sup_{w\in K}|P(w)|
\text{ for every }P\in\mathbb{C}[z_1,\ldots,z_n]
\right\}.$ The set $K$ is called polynomially convex if
$\widehat{K}=K$.Polynomial convexity plays an important role in several
complex variables, particularly in the theory of polynomial
approximation of continuous functions. However,
determining whether a given compact set is polynomially
convex is often a difficult problem. In this paper, we investigate the polynomial convexity of totally real closed discs contained in the boundaries of strongly pseudoconvex domains in $\mathbb{C}^2$. We construct an example of a non-polynomially convex totally real closed disc lying in the boundary of a bounded smoothly bounded strongly pseudoconvex domain. On the positive side, we establish sufficient conditions under which totally real discs are polynomially convex.

\medskip

Recall that a $C^1$-smooth real submanifold $M$ of
$\mathbb{C}^n$ is said to be \textit{totally real} at a point
$p\in M$ if $T_pM\cap iT_pM=\{0\},$ where $T_pM$ is regarded as a real linear subspace of
$\mathbb{C}^n$. The manifold $M$ is called \textit{totally real}
if it is totally real at every point. Totally real submanifolds play an important role in
polynomial and holomorphic approximation theory. Every
$C^1$-smooth totally real submanifold is locally
polynomially convex (see \cite{NirWell69}). Moreover,
every continuous function on a compact subset $K$ of
such a manifold can be uniformly approximated by
functions holomorphic on neighborhoods of $K$.
If $K$ is polynomially convex, then uniform
approximation by holomorphic polynomials is also
possible (see \cite{HoW68}).

By a \textit{totally real closed disc} in $\mathbb{C}^2$,
we mean a compact subset of a $C^1$-smooth
two-dimensional totally real submanifold
that is $C^1$-diffeomorphic to the closed
unit disc $\overline{\mathbb{D}}$. Despite the local polynomial convexity of totally real
submanifolds, a totally real closed disc in $\mathbb{C}^2$
need not be polynomially convex. The following classical
example of H\"ormander and Wermer \cite{HoW68} illustrates this
phenomenon.

\begin{example}\label{Exam: Wermer disc}
Let $f(z)=-(1+i)\overline z
+iz\overline z^{\,2}
+z^2\overline z^{\,3},$ and for $R>0,$ let us define
\begin{align*}
   K_{R}:=\{(z,f(z)):|z|\leq R\}. 
\end{align*}
It is easy to show that $K_{R}$ is a smooth
totally real closed disc. Since $f(z)=0$ whenever $|z|=1$, the polynomial hull
of $K_1$ contains the analytic disc $\{(z,0)\in\mathbb{C}^2:|z|\leq1\}.$ Thus $K_1$ is not polynomially convex. We will prove in Theorem~\ref{thm:main-truncated} that $K_R$
is polynomially convex for every $R\in(0,1)$. 
\end{example}

\medskip


Polynomial convexity of totally real discs has been
studied in several settings. For discs of graph type,
Duval \cite{JDuval88} established a sufficient
condition for polynomial convexity in terms of the
derivatives of the defining function. For totally
real discs contained in the boundaries of domains,
J\"oricke \cite{Jorc88} proved that every
$C^2$-smooth totally real closed disc in
$\partial\mathbb{B}^2$ is removable for CR-functions. Combining J\"oricke's removability theorem with a result
of Lupacciolu and Stout (see \cite{Sto93}), one obtains
that every $C^2$-smooth totally real closed disc contained
in $\partial\mathbb{B}^2$ is polynomially convex.
Stout also observed in \cite{Sto93} that a similar
argument applies to $C^2$-smooth totally real closed discs
contained in the boundary of a strictly pseudoconvex
domain $\Omega$, provided that $\overline{\Omega}$ is
polynomially convex. Later, Forstneri\v{c} and Stout \cite{ForSto1991}
extended this type of result to totally real discs
contained in the strongly pseudoconvex boundaries
of domains in two-dimensional Stein manifolds.

 \medskip

This naturally raises the question of whether the
polynomial convexity of $\overline{\Omega}$ is essential.
More precisely, we ask the following question:

\medskip

\noindent
\textbf{Question 1.}\label{Q: General case}
\textit{Is every $C^2$-smooth totally real closed disc
contained in the boundary of a smoothly bounded,
strongly pseudoconvex domain in $\mathbb{C}^2$
polynomially convex?}

\medskip

In our first result, we give a negative
answer to this question. More precisely, we prove the following:

\begin{theorem}\label{thm:counterexample}
There exist a bounded strongly pseudoconvex domain
$\Omega\subset\mathbb{C}^2$ with $C^\infty$-smooth boundary
and a real-analytic embedded totally real closed disc
$K\subset\partial\Omega$ such that $\widehat{K}\setminus\overline\Omega\neq\varnothing.$ In particular, neither $K$ nor $\overline\Omega$
is polynomially convex. Moreover, $\Omega$ is not Runge.
\end{theorem}

\medskip

Our construction is inspired by Wermer's classical
example \cite{Wermer1980}; see also the work of
Kim \cite{Kim1994}. We construct a
family of bounded domains with smooth strongly
pseudoconvex boundaries containing a non-polynomially
convex totally real closed disc $K_1$ (see, Example \ref{Exam: Wermer disc}). Moreover, the
polynomial hull of this disc extends beyond the
closures of these domains.

\medskip

Our second result concerns polynomial convexity of totally real
discs contained in Levi-flat hypersurfaces. 
For results on the polynomial convexity of totally real discs
contained in Levi-flat hypersurfaces, see
\cite{Alex91,GorMon2025}.
 In \cite{GorMon2025}, polynomial convexity of totally
real discs contained in several classes of Levi-flat
hypersurfaces in $\mathbb{C}^2$ was established.
The following conjecture was also proposed
in \cite[Conjecture 5.5]{GorMon2025}:
\textit{Any totally real disc lying in a nonsingular
Levi-flat hypersurface in $\mathbb{C}^2$ is
polynomially convex.}

\smallskip

\noindent We show that the conjecture fails when Levi-flat
hypersurfaces are allowed to be defined only on
neighborhoods of the discs. We construct
a nonsingular real-analytic Levi-flat hypersurface
in an open subset of $\mathbb{C}^2$ containing a
smooth totally real closed disc that is not
polynomially convex. More precisely, using H\"ormander--Wermer disc $K_1$, we prove the following result.

\medskip
\begin{theorem}\label{T: Disc_Levi-flat}
There exist non-Runge domain $V\subset\mathbb{C}^2$,
a nonsingular real-analytic Levi-flat hypersurface
$M\subset V$, and a real-analytic embedded
totally real closed disc $K\subset M$ such that $K$ is not polynomially convex.
\end{theorem}

\begin{remark}
 Theorem~\ref{T: Disc_Levi-flat} does not settle the conjecture
for globally closed nonsingular
Levi-flat hypersurfaces in
$\mathbb{C}^2$.   
\end{remark}

\medskip

We next study the polynomial convexity of certain
totally real discs of graph type. Although the
classical H\"ormander--Wermer disc $K_1$ is not
polynomially convex, it is locally polynomially
convex, since it is a smooth totally real disc.
In this direction, Gorai \cite[Example~5.1]{SG17}
showed that $K_R$ is polynomially convex for every
$R\in(0,1/\sqrt{3}]$.

\medskip

\noindent Another example of a non-polynomially convex totally
real disc was given by Forstneri\v{c}
\cite{Forstneric1992}: for $R\in (0,1]$, let us define
\begin{align*}
E_R:=\left\{
(z,g(z))\in\mathbb{C}^2:
g(z)=\overline z(1-|z|^2)e^{i|z|^2},
\ |z|\leq R
\right\}.
\end{align*}
Since $g(z)=0$ whenever $|z|=1$, the polynomial hull
of $E_1$ contains the analytic disc $\{(z,0)\in\mathbb{C}^2:|z|\leq1\}.$ Thus $E_1$ is not polynomially convex. On the other hand,
Gorai \cite[Remark 5.2]{SG17} showed that $E_R$ is polynomially
convex for $0<R<1/16\sqrt{11}$. In the following theorem, we improve these results
by proving that both $K_R$ and $E_R$ are polynomially
convex for every $R<1$. Our proof uses a different
approach based on Stolzenberg's fiber criterion.

\begin{theorem}\label{thm:main-truncated}
For every $R\in (0,1)$, the totally real discs
$K_R$ and $E_R$ are polynomially convex.
\end{theorem}

\medskip

Finally, we return to the question
of finding sufficient conditions
for polynomial convexity of totally
real discs contained in strongly
pseudoconvex boundaries. As mentioned earlier, the result
of Stout \cite{Sto93}, Forstneri\v{c} and Stout
\cite{ForSto1991} gives
holomorphic convexity of such discs
relative to the closure of the
ambient domain. The difficulty
is that the closure of a strongly
pseudoconvex domain need not be
polynomially convex (see, \cite[Example 2.7]{AJI1993}).

\smallskip

In \cite{ChatterjeeGorai2025}, Chatterjee and Gorai
studied polynomial convexity of the closures of
strongly pseudoconvex domains invariant under the
positive-time flows of certain holomorphic vector fields. Motivated by their work, we
consider a condition on a holomorphic vector field
along a totally real disc in the boundary, without
requiring invariance of the entire domain. We prove the following result:

\begin{theorem}\label{thm:main-inward}
Let $\Omega\Subset\mathbb{C}^2$ be a bounded Runge
domain with $C^2$-smooth strongly
pseudoconvex boundary, and let
$K\subset\partial\Omega$ be a compact
$C^2$-smooth embedded totally real
closed disc. Suppose there exists a holomorphic vector
field $X$ on a neighborhood of $\widehat{K}$
whose local holomorphic flow
$\{\phi_t\}_{|t|<\varepsilon}$ is defined
in a neighborhood of $\widehat{K}$
and satisfies
\begin{align}\label{eq:inward-flow-condition}
\phi_t(K)\subset\Omega
\qquad\text{for every }0<t<\varepsilon.
\end{align}
Then $K$ is polynomially convex.
\end{theorem}

\begin{remark}
    In view of Result~\ref{res:local-holomorphic-flow} such a holomorphic flow $\{\phi_t\}_{|t|<\varepsilon}$ defined in a neighborhood of $\widehat{K}$ exists.
\end{remark}

The following result establishes the converse of
Theorem~\ref{thm:main-inward} and provides several
equivalent conditions for the polynomial convexity
of totally real closed discs.

\begin{theorem}\label{T:vector-field-equivalence}
Let $\Omega\Subset\C^2$ be a bounded Runge
domain with $C^\infty$-smooth strongly
pseudoconvex boundary, and let
$K\subset\partial\Omega$ be a $C^2$-smooth
embedded totally real closed disc.
Let $\rho$ be any smooth defining function
for $\Omega$. Then the following statements are equivalent:
\begin{enumerate}
\item[(i)] $K$ is polynomially convex.

\item[(ii)] There exists an entire holomorphic
vector field $X:\mathbb{C}^2\to\mathbb{C}^2$ such that $d\rho_p(X(p))<0$
for every $p\in K.$

\item[(iii)] There exists a polynomial
holomorphic vector field $X:\mathbb{C}^2\to\mathbb{C}^2$
such that $d\rho_p(X(p))<0$
for every $p\in K.$

\item[(iv)] There exists a local holomorphic flow
$\{\phi_t\}_{|t|<\varepsilon}$, generated by a holomorphic
vector field defined on a neighborhood of $\widehat{K}$,
such that $\phi_t(K)\subset\Omega$
for every $t\in(0,\varepsilon).$

\item[(v)] There exists a polynomial
holomorphic vector field $X$ whose local
holomorphic flow $\{\phi_t\}$ satisfies $\phi_t(K)\subset\Omega$ for all sufficiently small $t>0$.
\end{enumerate}
\end{theorem}

\section{Preliminaries}\label{S:technical}

\subsection{Holomorphic vector fields and local flows}

A holomorphic vector field $V$ on a domain
$U\subset\mathbb{C}^n$ is a real vector field of the form
\begin{align*}
V=\sum_{j=1}^n\left(
a_j(z)\frac{\partial}{\partial x_j}
+b_j(z)\frac{\partial}{\partial y_j}
\right),
\end{align*}
where $z_j=x_j+iy_j$ and
$f_j:=a_j+ib_j$ is holomorphic on $U$
for every $j=1,\ldots,n$. Equivalently, $V$ is determined by the
holomorphic mapping
$f=(f_1,\ldots,f_n):U\to\mathbb{C}^n$.
The local holomorphic flow generated by $V$
is the solution of
\begin{align*}
\frac{d}{dt}\phi_t(z)=f(\phi_t(z)),
\qquad
\phi_0(z)=z,
\end{align*}
for all sufficiently small complex $t$ and
all $z$ in a suitable neighborhood.

By the local existence and uniqueness theorem
for holomorphic ordinary differential equations,
every holomorphic vector field admits a unique
local holomorphic flow. We have the following result from Bracci \cite[Theorem~1.21]{Bracci2010}.
\begin{result}
\label{res:local-holomorphic-flow}
Let $V\subset\mathbb{C}^n$ be an open set, and let
$X:V\to\mathbb{C}^n$ be a holomorphic vector field.
For every compact set $L\subset V$, there exist
$\varepsilon>0$, an open neighborhood $U$ of $L$,
and a unique real-analytic mapping $\phi:(-\varepsilon,\varepsilon)\times U
\longrightarrow V$
such that, for each fixed
$t\in(-\varepsilon,\varepsilon)$,
the mapping $\phi_t:=\phi(t,\cdot)$ is
holomorphic on $U$ and satisfies
\begin{align*}
\begin{cases}
\displaystyle
\frac{d}{dt}\phi_t(z)=X(\phi_t(z)),\\[1pt]
\phi_0(z)=z,\quad z\in U.
\end{cases}
\end{align*}
\end{result}

Let $\Omega\subset\mathbb{C}^n$ be an open set and
$E\subset\Omega$ be compact. We denote by
$\widehat{E}_{\mathcal{O}(\Omega)}$ the
$\mathcal{O}(\Omega)$-hull of $E$, defined by
\begin{align*}
\widehat{E}_{\mathcal{O}(\Omega)}
:=
\left\{
z\in\Omega:
|f(z)|\leq\sup_{w\in E}|f(w)|
\text{ for every }f\in\mathcal{O}(\Omega)
\right\}.
\end{align*}
We say that $E$ is $\mathcal{O}(\Omega)$-convex
if $\widehat{E}_{\mathcal{O}(\Omega)}=E.$

\begin{remark}
When $\Omega=\mathbb{C}^n$, the $\mathcal{O}(\Omega)$-hull
of $E$ coincides with its polynomial convex hull.
Indeed, every entire function can be approximated
uniformly on compact subsets of $\mathbb{C}^n$ by
holomorphic polynomials. Thus, $\widehat E_{\mathcal{O}(\mathbb{C}^n)}=\widehat E.$
\end{remark}

Let $K$ be a compact subset of $\mathbb{C}^n$. We denote by
$\mathcal{C}(K)$ the space of all complex-valued continuous functions
on $K$, equipped with the supremum norm
$\|f\|_K:=\sup_{z\in K}|f(z)|$.
Let $\mathcal{P}(K)$ denote the uniform closure in $\mathcal{C}(K)$ of
the restrictions to $K$ of holomorphic polynomials.

The following approximation theorem is due to
O'Farrell, Preskenis, and Walsh \cite{OPW1983}.

\begin{result}[O'Farrell--Preskenis--Walsh]
\label{OPW1983}
Let $K$ be a compact polynomially convex subset of
$\mathbb{C}^n$, and let $E$ be a closed subset of $K$.
Suppose that $K\setminus E$ is locally contained
in a $C^1$-smooth totally real submanifold of $\mathbb{C}^n$.
Then
\begin{align*}
\mathcal{P}(K)=\{f\in \mathcal{C}(K):f|_E\in \mathcal{P}(E)\}.
\end{align*}
\end{result}

Recall that a domain $\Omega\subset\mathbb{C}^2$ is \emph{Runge} if
holomorphic functions on $\Omega$ can be approximated uniformly
on compact subsets by holomorphic polynomials
in $\mathbb{C}^2$. 

\smallskip

The following result is from Gorai--Mondal \cite[Theorem~1.6]{GorMon2025}

\begin{result}
\label{R:Levi-flat-polynomial-convexity}
Let $\Omega\subset\mathbb{C}^2$ be a Runge domain,
and let $h\in\mathcal O(\Omega)$. Suppose that $M:=\{z\in\Omega:|h(z)|=1\}$
satisfies $dh(z)\neq0$ for every $z\in M$. Then every $C^2$-smooth totally real closed disc contained in $M$ is polynomially convex.
\end{result}

\medskip

We shall use the following result of
Forstneri\v{c} and Stout \cite[Corollary II]{ForSto1991}.

\begin{result}\label{R:Holomorphic Convexity}
   Let $D$ be a strongly pseudoconvex domain in a
two-dimensional Stein manifold $\mathcal{M}$.
Then every two-dimensional compact totally real disc
$\Delta$ of class $C^1$ contained in $\partial D$
is removable and hence
$\mathcal{O}(\overline{D})$-convex. 
\end{result}

\noindent The following is due to Stolzenberg (see \cite[Lemma 2.3]{Gorai2014}).
	\begin{result}[Stolzenberg]\label{R:PolyCnvx_Stolzenberg}
		Let $E$ be a compact subset of $\mathbb{C}^{n}.$ Assume that $\poly(E)$ contains a function $\phi$ such that $\phi(E)$ has empty interior and $\mathbb{C}\setminus\phi(E)$ is connected. Then $E$ is polynomially convex if and only if $\phi^{-1}(c)\cap E$ is polynomially convex for each $c\in \phi(E).$ 	
	\end{result}

We also make use of the following two elementary observations.

\begin{lemma}\label{lem:holomorphic-image-hull}
Let $E\subset\mathbb{C}^n$ be compact, and let
$F:U\to\mathbb{C}^m$ be a holomorphic mapping,
where $U$ is an open neighborhood of
$\widehat E$. Then $F(\widehat E)\subset\widehat{F(E)}.$
\end{lemma}

\begin{proof}
Let $q\in\widehat E$ and let $P$ be an
arbitrary holomorphic polynomial on $\mathbb{C}^m$.
Since $P\circ F$ is holomorphic on a
neighborhood of $\widehat E$, and
$\widehat E$ is polynomially convex,
the Oka--Weil theorem gives a sequence
of holomorphic polynomials $Q_j$ on
$\mathbb{C}^n$ such that$Q_j\longrightarrow P\circ F
\quad\text{uniformly on }\widehat E.$ Since $q\in\widehat E$, we have
\begin{align*}
|Q_j(q)|\leq\sup_{p\in E}|Q_j(p)|.
\end{align*}
Letting $j\to\infty$, we obtain $|P(F(q))|
\leq\sup_{p\in E}|P(F(p))|
=\sup_{w\in F(E)}|P(w)|.$ Since $P$ was arbitrary, it follows that
$F(q)\in\widehat{F(E)}$. Hence $F(\widehat E)\subset\widehat{F(E)}.$
\end{proof}

\begin{lemma}
\label{lem:hull-containment}
Let $\Omega\Subset\mathbb{C}^2$ be a bounded
domain with $C^\infty$-smooth strongly
pseudoconvex boundary, and let
$K\subset\partial\Omega$ be a $C^2$-smooth
embedded totally real closed disc. If $\widehat{K}\subset\overline\Omega,$ then $K$ is polynomially convex.
\end{lemma}

\begin{proof}
By Result~\ref{R:Holomorphic Convexity}, the disc $K$
is convex with respect to the algebra
$\mathcal O(\overline\Omega)$ of functions
holomorphic on neighborhoods of
$\overline\Omega$. Thus, for every
$q\in\overline\Omega\setminus K$,
there exists a function $g$, holomorphic
on a neighborhood of $\overline\Omega$,
such that
\begin{equation}\label{eq:FS-separation}
|g(q)|>\sup_{p\in K}|g(p)|.
\end{equation}

Suppose, to the contrary, that $K$ is not polynomially convex. Choose $q\in\widehat{K}\setminus K.$ By hypothesis, $q\in\overline\Omega\setminus K.$ Hence there exists a function $g$ satisfying
\eqref{eq:FS-separation}.

\medskip

\noindent On the other hand, $\widehat{K}$ is a polynomially convex compact set, and $\widehat{K}\subset\overline\Omega.$ Therefore $g$ is holomorphic on a
neighborhood of $\widehat{K}$. By the Oka--Weil approximation theorem, there exists a sequence of holomorphic polynomials $\{P_j\}$ converging uniformly
to $g$ on $\widehat{K}$. Since $q\in\widehat{K}$, we have $|P_j(q)|\leq\sup_K|P_j|$ for every $j.$ Hence $|g(q)|\leq\sup_K|g|,$ contradicting \eqref{eq:FS-separation}. Therefore $K$ is polynomially convex.
\end{proof}

\medskip


Recall that $K_1=\{(z,f(z)):|z|\leq1\},$ where $f(z)=-(1+i)\overline z
+iz\overline z^{\,2}
+z^2\overline z^{\,3}.$ For simplicity, we write $K:=K_1$ throughout
the remainder of this section.

\smallskip

\noindent Consider the holomorphic polynomial $F(z,\zeta)=-(1+i)\zeta+iz\zeta^2+z^2\zeta^3,$ so that $f(z)=F(z,\overline z)$.
Following Wermer \cite{Wermer1980}
(see also \cite{Kim1994}), define the
holomorphic polynomial mapping
\begin{align*}
\Phi:\mathbb{C}^2\longrightarrow\mathbb{C}^2,
\qquad
\Phi(u,v)
=\bigl(u+iv,F(u+iv,u-iv)\bigr).
\end{align*}
Let
\begin{equation}\label{E:Real Disc}
\Delta_{\R}
:=\{(u,v)\in\R^2:u^2+v^2\leq1\}.\quad \text{ Consequently }, \Phi(\Delta_{\R})=K.
\end{equation}

\smallskip

We next show that $\Phi$ is locally
biholomorphic at every point of
$\Delta_{\R}$. Clearly, $\Phi$ is
injective on $\Delta_{\R}$, since
its first component $u+iv$ uniquely
determines $(u,v)\in\R^2$. By the chain rule, $\det D\Phi(u,v)
=-2iF_\zeta(u+iv,u-iv),$ where $F_\zeta(z,\zeta)
=-(1+i)+2iz\zeta+3z^2\zeta^2.$ For $(u,v)\in\Delta_{\R}$, putting
$z=u+iv$, we obtain
\begin{align*}
F_\zeta(z,\overline z)
=-(1+i)+2i|z|^2+3|z|^4
=(3|z|^4-1)+i(2|z|^2-1)
\neq0.
\end{align*}
Consequently,
\begin{equation}\label{eq:phi-jacobian}
\det D\Phi(u,v)\neq0
\qquad\text{for every }(u,v)\in\Delta_{\R}.
\end{equation}
Hence, by the holomorphic inverse function
theorem, $\Phi$ is locally biholomorphic
at every point of $\Delta_{\R}$.

\medskip

With the above notation, we now establish a few lemmas
that will be used in the proof of
Theorem~\ref{thm:counterexample}.

\begin{lemma}\label{L:Existence_Bihoplo}
There exists an open
neighborhood $U$ (in $\mathbb{C}^2$) of $\Delta_{\R}$ such that $\Phi|_U:U\longrightarrow\Phi(U)
\quad\text{is biholomorphic}.$
\end{lemma}

\begin{proof}
  For $\varepsilon>0$, define
\begin{align*}
U_\varepsilon
=
\left\{
p\in\mathbb{C}^2:
\operatorname{dist}(p,\Delta_{\R})
<\varepsilon
\right\}.
\end{align*}

We first claim that $\Phi$ is injective on
$U_\varepsilon$ for some sufficiently small
$\varepsilon>0$. Suppose, to the contrary, that no such
$\varepsilon$ exists. Then, for every
positive integer $j$, there exist distinct
points $a_j,b_j\in U_{1/j}$ such that
\begin{align}\label{eq:phi-equal-points}
a_j\neq b_j,
\qquad
\Phi(a_j)=\Phi(b_j).
\end{align}
\noindent By the definition of $U_{1/j}$, we can choose
$p_j,q_j\in\Delta_{\R}$ such that $\|a_j-p_j\|<\frac1j,
\qquad
\|b_j-q_j\|<\frac1j.$ Since $\Delta_{\R}$ is compact, after passing
to subsequences, we may assume that $p_j\longrightarrow p\in\Delta_{\R},
\qquad
q_j\longrightarrow q\in\Delta_{\R}.$ It follows that $a_j\longrightarrow p,
\qquad
b_j\longrightarrow q.$ Since $\Phi$ is continuous and
$\Phi(a_j)=\Phi(b_j)$, passing to the limit
gives $\Phi(p)=\Phi(q).$ But $\Phi$ is injective on $\Delta_{\R}$.
Therefore, $p=q.$ Consequently,
\begin{equation}\label{eq:common-limit}
a_j\longrightarrow p,
\qquad
b_j\longrightarrow p.
\end{equation}

On the other hand, by
\eqref{eq:phi-jacobian}, we have $\det D\Phi(p)\neq0.$ The holomorphic inverse function theorem
therefore provides an open neighborhood
$V_p$ of $p$ such that $\Phi|_{V_p}$ is biholomorphic onto its image. In particular,
$\Phi$ is injective on $V_p$. By \eqref{eq:common-limit}, for all sufficiently
large $j$, both $a_j$ and $b_j$ belong to $V_p$.
Since $\Phi(a_j)=\Phi(b_j),$ injectivity of $\Phi|_{V_p}$ implies $a_j=b_j,$ contradicting \eqref{eq:phi-equal-points}. This contradiction proves that there exists
$\varepsilon_1>0$ such that $\Phi|_{U_{\varepsilon_1}}
\quad\text{is injective}.$

\medskip

\noindent Next, since $\det D\Phi$ is continuous and
does not vanish on the compact set
$\Delta_{\R}$, there exists
$\varepsilon_2>0$ such that $\det D\Phi(p)\neq0$ for every $p\in U_{\varepsilon_2}.$ Choose $0<\varepsilon<
\min\{\varepsilon_1,\varepsilon_2\}$ and set $U=U_\varepsilon.$ Then $\Phi$ is injective on $U$, and $\det D\Phi(p)\neq0
\qquad
\text{for every }p\in U.$ By the holomorphic inverse function theorem,
$\Phi$ is locally biholomorphic at every
point of $U$. In particular, $\Phi(U)$
is open. Thus $\Phi|_U:U\longrightarrow\Phi(U)
\text{ is biholomorphic}.$
\end{proof}

\medskip

\begin{lemma}\label{L:Construction_D_c}
There exists a family $\{D_c\}_{c>0}$ of smoothly bounded, convex, strongly pseudoconvex
domains in $\mathbb{C}^2$ such that $\Delta_{\R}\subset\partial D_c$
for every $c>0.$
\end{lemma}

\begin{proof}

We recall the following function from \cite{Chi2016}: let $\chi:\R\to[0,\infty)$ be defined by
\begin{align*}
\chi(t)=
\begin{cases}
0, & t\leq0,\\[4pt]
\displaystyle\int_0^t e^{-1/s}\,ds, & t>0.
\end{cases}
\end{align*}
Then $\chi$ is $C^\infty$-smooth, convex,
and nondecreasing on $\R$. Moreover,
$\chi(t)>0$ for every $t>0$, and
$\chi(t)\to+\infty$ as $t\to+\infty$. Indeed, for $t>0$,
\begin{align*}
\chi'(t)=e^{-1/t}>0,
\qquad
\chi''(t)=\frac{e^{-1/t}}{t^2}>0.
\end{align*}


\medskip

Write $u=x_1+iy_1,\quad
v=x_2+iy_2,$ where $x_1,x_2,y_1,y_2\in\R$. For every $c>0$, define
\begin{align*}
\rho_c(u,v)
:=
y_2+c(y_1^2+y_2^2)
+\chi(x_1^2+x_2^2-1)
\end{align*}
and
\begin{align*}
D_c
:=
\{(u,v)\in\mathbb{C}^2:\rho_c(u,v)<0\}.
\end{align*}

It is easy to verify that $D_c$ is nonempty and bounded.
Moreover, since $\rho_c$ is continuous, $D_c$ is open.

\medskip

Note that $\rho_c$ is a convex function on $\R^4$, and hence $D_c$ is convex. Indeed, put $x=(x_1,x_2)\in\R^2,
\qquad
y=(y_1,y_2)\in\R^2.$ Then $\rho_c(x,y)
=
y_2+c|y|^2+\chi(|x|^2-1).$ The function $y\mapsto y_2$ is affine
and hence convex, while
$y\mapsto c|y|^2$ is convex because
$c>0$. 
Since $\chi$ is convex and nondecreasing, and
$x\mapsto |x|^2-1$ is convex, the composition
$x\mapsto\chi(|x|^2-1)$ is convex. Therefore, $\rho_c$ is convex.

\medskip

\noindent
\textit{\bf Claim 1:} $\rho_c$ is Strictly plurisubharmonic. We next prove that $\rho_c$ is
strictly plurisubharmonic on $\mathbb{C}^2$. Recall that for a real-valued
$C^2$-function $\rho$ on $\mathbb{C}^2$,
its Levi form is
\begin{align*}
\mathcal L_\rho(p;\xi)
=
\sum_{j,k=1}^{2}
\frac{\partial^2\rho}
{\partial z_j\partial\overline z_k}(p)
\,\xi_j\overline{\xi_k},
\end{align*}
where $(z_1,z_2)=(u,v),
\qquad
\xi=(\xi_1,\xi_2)\in\mathbb{C}^2.$

All mixed derivatives of $c(y_1^2+y_2^2)$ involving
different coordinates vanish.
Consequently,
\begin{align}\label{E:Levi_2nd fun}
   \mathcal L_{c(y_1^2+y_2^2)}(p;\xi)
=
\frac c2\left(
|\xi_1|^2+|\xi_2|^2
\right). 
\end{align}
\noindent By a simple computation, we get the Levi form of $h(u,v):=\chi(x_1^2+x_2^2-1)$ is
\begin{align}\label{E:Levi_3rd fun}
\notag \mathcal L_h(p;\xi)
&=
\sum_{j,k=1}^{2}
\left(
\chi''(s)x_jx_k
+\frac12\chi'(s)\delta_{jk}
\right)
\xi_j\overline{\xi_k}\\ 
&=
\chi''(s)
\left|x_1\xi_1+x_2\xi_2\right|^2
+
\frac{\chi'(s)}2
\left(|\xi_1|^2+|\xi_2|^2\right).
\end{align}
Since $\chi'(s)\geq0$ and
$\chi''(s)\geq0$, we have $\mathcal L_h(p;\xi)\geq0.$ Combining (\ref{E:Levi_2nd fun}) and (\ref{E:Levi_3rd fun}),
we obtain $\mathcal L_{\rho_c}(p;\xi)
=
\frac c2
\left(|\xi_1|^2+|\xi_2|^2\right)
\nonumber+
\frac{\chi'(s)}2
\left(|\xi_1|^2+|\xi_2|^2\right)
\nonumber+
\chi''(s)
\left|x_1\xi_1+x_2\xi_2\right|^2.$ In particular, $\mathcal L_{\rho_c}(p;\xi)
\geq
\frac c2\|\xi\|^2.$
Since $c>0$, the right-hand side
is strictly positive whenever
$\xi\neq0$. Therefore $\rho_c$ is strictly
plurisubharmonic on $\mathbb{C}^2.$

\noindent This establishes the claim 2.\hfill $\blacktriangleleft$

\medskip

\medskip

\noindent
\textit{\bf Claim 2.} The real disc $\Delta_{\R}$ lies in $\partial D_c$. Recall that $\Delta_{\R}
=
\left\{
(u,v)\in\R^2:
u^2+v^2\leq1
\right\}.$ Let $(u,v)\in\Delta_{\R}$.
Since $u,v$ are real, we have $y_1=y_2=0.$ Moreover, $x_1^2+x_2^2\leq1.$ By the definition of $\chi$, $\chi(x_1^2+x_2^2-1)=0.$ Therefore $\rho_c(u,v)
=0+c(0)+0=0.$ Hence, $ \Delta_{\R}\subset\partial D_c$
for every $c>0.$ 

\noindent This proves the { Claim 2}.$\hfill{\blacktriangleleft}$

\smallskip

 Since $\chi$ is $C^\infty$,
the function $\rho_c$ is
$C^\infty$ on $\mathbb{C}^2$. A direct differentiation gives us that $\nabla_{\R}\rho_c(p)\neq0$
whenever $\rho_c(p)=0.$ Therefore,
$\{\rho_c=0\}$ is a
$C^\infty$-smooth real hypersurface, and  $\partial D_c=\{\rho_c=0\},$ $\overline{D_c}=\{\rho_c\leq0\}.$ 

\noindent Thus, { Claim 1}, and {Claim 2} establishes the result.

\end{proof}

\begin{lemma}\label{L:Hausdorff Convergence}
As $c\to+\infty$, the compact sets
$\overline{D}_c$ converge to $\Delta_{\R}$
in the Hausdorff metric.
\end{lemma}

\begin{proof}
Observe that for $c_2>c_1$, $\rho_{c_2}(u,v)-\rho_{c_1}(u,v)
=(c_2-c_1)(y_1^2+y_2^2)\geq0.$ Hence $\overline{D}_{c_2}\subset\overline{D}_{c_1}.$

\smallskip

\noindent {\bf Claim.} We claim that $\bigcap_{c\geq1}\overline{D}_c=\Delta_{\R}.$

Indeed, if
$p=(x_1+iy_1,x_2+iy_2)$ belongs to
this intersection, then $y_2+c(y_1^2+y_2^2)
+\chi(x_1^2+x_2^2-1)\leq0$ for every $c\geq1$.
Letting $c\to+\infty$ gives
$y_1=y_2=0$. Consequently, $\chi(x_1^2+x_2^2-1)=0,$ and therefore $x_1^2+x_2^2\leq1$.
The reverse inclusion follows from
$\Delta_{\R}\subset\overline{D}_c$.

\noindent This proves the claim. $\hfill{\blacktriangleleft}$

\smallskip

Recall that for nonempty compact
sets $A,B\subset\mathbb{C}^2$,
\begin{align*}
d_H(A,B)
=
\max\left\{
\sup_{a\in A}\operatorname{dist}(a,B),
\sup_{b\in B}\operatorname{dist}(b,A)
\right\}.
\end{align*}

\noindent Since $\Delta_{\R}\subset\overline{D}_c$ for every $c>0$, we have $\sup_{p\in\Delta_{\R}}
\operatorname{dist}(p,\overline{D}_c)=0.$ Thus it suffices to prove that $\sup_{p\in\overline{D}_c}
\operatorname{dist}(p,\Delta_{\R})
\longrightarrow 0.$ Suppose now that the asserted Hausdorff
convergence fails. Then,
there exist $\varepsilon>0$,
a sequence $c_j\to+\infty$, and
points $p_j\in\overline{D}_{c_j}$
such that
\begin{align*}
\operatorname{dist}(p_j,\Delta_{\R})
\geq\varepsilon.
\end{align*}
For sufficiently large $j$, we have
$p_j\in\overline{D}_1$.
By compactness of $\overline{D}_1$,
after passing to a subsequence,
we may assume $p_j\to p$. For every fixed $c\geq1$, we have
$p_j\in\overline{D}_c$ for all
sufficiently large $j$. Hence
$p\in\overline{D}_c$ for every $c\geq1$.
Thus $p\in\Delta_{\R}$ (by above claim), contradicting $\operatorname{dist}(p_j,\Delta_{\R})
\geq\varepsilon.$ Therefore, $d_H(\overline{D}_c,\Delta_{\R})
\longrightarrow0$ as $c\to+\infty.$
\end{proof}

\begin{lemma}\label{L:Construction_Omega_c}
For all sufficiently large $c>0$, the domain
$\Omega_c:=\Phi(D_c)$ is bounded and has
$C^\infty$-smooth strongly pseudoconvex boundary.
Moreover, $K\subset\partial\Omega_c.$
\end{lemma}

\begin{proof}
Recall from Lemma~\ref{L:Existence_Bihoplo} that there exists an
open neighborhood $U$ of $\Delta_{\R}$
such that $\Phi|_U:U\longrightarrow\Phi(U)$ is biholomorphic.

\smallskip

\noindent{\bf Claim.} $\overline{D}_c\subset U$ for all sufficiently large $c$.

Indeed, since $\Delta_{\R}$ is compact
and $U$ is an open neighborhood of
$\Delta_{\R}$, there exists
$\varepsilon_0>0$ such that $\left\{
p\in\mathbb{C}^2:
\operatorname{dist}(p,\Delta_{\R})
<\varepsilon_0
\right\}
\subset U.$ By Lemma~\ref{L:Hausdorff Convergence}, there exists
$c_0>0$ such that, whenever $c\geq c_0$, $\sup_{p\in\overline{D_c}}
\operatorname{dist}(p,\Delta_{\R})
<\varepsilon_0.$ Therefore,
\begin{align*}
\overline{D}_c\subset U
\qquad\text{for every }c\geq c_0.
\end{align*}
\noindent This proves the claim.$\hfill{\blacktriangleleft}$

\smallskip

For such $c$, define $\Omega_c=\Phi(D_c).$ Since $\Phi|_U$ is biholomorphic
and $\overline{D}_c\subset U$,
the restriction $\Phi|_{D_c}:D_c\longrightarrow\Omega_c$ is biholomorphic. Moreover, $\partial\Omega_c=\Phi(\partial D_c).$

\smallskip

\noindent We observe that the function $\widetilde\rho_c
=
\rho_c\circ(\Phi|_U)^{-1}$ is a $C^\infty$-smooth defining function for $\Omega_c$ in a
neighborhood of its boundary. If $w=\Phi(p)$ and
$\eta\in\mathbb{C}^2\setminus\{0\}$,
the holomorphic chain rule gives $\mathcal L_{\widetilde\rho_c}(w;\eta)
=
\mathcal L_{\rho_c}
\left(
p;
D\Phi^{-1}(w)\eta
\right).$ Since $D\Phi^{-1}(w)$ is invertible
and $\rho_c$ is strictly
plurisubharmonic, we have $\mathcal L_{\widetilde\rho_c}(w;\eta)>0.$ Also, since $\Phi$ is a
diffeomorphism on $U$, the gradient
of $\widetilde\rho_c$ does not vanish
on $\partial\Omega_c$. Thus $\Omega_c$ has
$C^\infty$-smooth, strongly
pseudoconvex boundary.

\medskip

\noindent 
Furthermore, $\Omega_c$ is bounded
because $\overline{D}_c$ is compact
and $\Phi$ is continuous on $(\overline{D}_c\subset )U$. Finally, using
Lemma~\ref{L:Construction_D_c} and
$\Phi(\Delta_{\R})=K$, we obtain $K=\Phi(\Delta_{\R})
\subset\Phi(\partial D_c)
=\partial\Omega_c.$
This completes the construction.
\end{proof}

\section{Proof of Theorem~\ref{thm:counterexample}}
First, we establish the following: 
\begin{align}\label{E:Lipschitz-estimate}
d_H(\overline{\Omega}_c,K)
\leq M\varepsilon_c
\end{align}
for some constant $M>0$ independent of
all sufficiently large $c$, where $\varepsilon_c
:=d_H(\overline{D}_c,\Delta_{\R}).$  Recall that $\Phi$ is biholomorphic on an
open neighborhood $U$ of $\Delta_{\R}$.
Since $\Delta_{\R}$ is compact, there
exists $r>0$ such that
\begin{align*}
T_r:=
\left\{
p\in\mathbb{C}^2:
\operatorname{dist}(p,\Delta_{\R})\leq r
\right\}
\subset U.
\end{align*}
Observe that $T_r$ is compact and convex. By Lemma~\ref{L:Hausdorff Convergence}, $\varepsilon_c\longrightarrow0
\quad\text{as }c\to+\infty.$ Therefore, there exists $c_0>0$ such that
$\varepsilon_c<r
~\text{for all }c\geq c_0.$ In particular, $\overline{D}_c\subset T_r
\qquad(c\geq c_0).$ Since $\Phi$ is holomorphic on $U$,
its derivative $D\Phi$ is continuous.
By compactness of $T_r$, there exists
a constant $M>0$ such that
\begin{align*}
\sup_{p\in T_r}
\|D\Phi(p)\|_{\mathrm{op}}\leq M.
\end{align*}
Here $\|\cdot\|_{\mathrm{op}}$ denotes
the operator norm induced by the Euclidean
norm on $\mathbb{C}^2$. We note that $T_r$ is fixed,
independently of $c$. Consequently,
the constant $M$ is also independent
of $c$. Thus $\Phi$ is Lipschitz continuous
on $T_r$ with Lipschitz constant $M$:
\begin{equation}\label{eq:step4-Phi-Lipschitz}
\|\Phi(a)-\Phi(b)\|
\leq M\|a-b\|
\qquad a,b\in T_r.
\end{equation}

We now estimate the Hausdorff distance
between $\overline{\Omega}_c$ and $K$. Let $w\in\overline{\Omega}_c$ be
arbitrary. Then there exists
$p\in\overline{D}_c$ such that $w=\Phi(p).$ Since $\Delta_{\R}$ is compact,
there exists $p_0\in\Delta_{\R}$
satisfying $\|p-p_0\|
=
\operatorname{dist}(p,\Delta_{\R}).$ By the definition of the Hausdorff
distance,
\begin{align*}
\operatorname{dist}(p,\Delta_{\R})
\leq
d_H(\overline{D}_c,\Delta_{\R})
=\varepsilon_c.
\end{align*}
Hence $\|p-p_0\|\leq\varepsilon_c.$ Since $p\in\overline{D}_c\subset T_r
\quad\text{and}\quad
p_0\in\Delta_{\R}\subset T_r,$ we may apply
\eqref{eq:step4-Phi-Lipschitz} to obtain $\|\Phi(p)-\Phi(p_0)\|
\leq M\|p-p_0\|
\leq M\varepsilon_c.$ But $\Phi(p_0)\in\Phi(\Delta_{\R})=K.$ Consequently, $\operatorname{dist}(w,K)
\leq
\|\Phi(p)-\Phi(p_0)\|
\leq M\varepsilon_c.$ Since $w\in\overline{\Omega}_c$
was arbitrary, we conclude that
\begin{equation}\label{eq:step4-one-sided-Hausdorff}
\sup_{w\in\overline{\Omega}_c}
\operatorname{dist}(w,K)
\leq M\varepsilon_c.
\end{equation}

On the other hand, we have already
proved that $\Delta_{\R}\subset\partial D_c.$ Applying $\Phi$, we obtain $K=\Phi(\Delta_{\R})
\subset\Phi(\partial D_c)
=\partial\Omega_c
\subset\overline{\Omega}_c.$ Therefore,
\begin{align}\label{eq:step4-reverse-Hausdorff}
\sup_{z\in K}
\operatorname{dist}(z,\overline{\Omega}_c)
=0.
\end{align}
\noindent Combining
(\ref{eq:step4-one-sided-Hausdorff})
and
(\ref{eq:step4-reverse-Hausdorff}),
we obtain
\begin{align}\label{eq:step4-one-sided}
d_H(\overline{\Omega}_c,K)
=
\max\left\{
\sup_{w\in\overline{\Omega}_c}
\operatorname{dist}(w,K),
\sup_{z\in K}
\operatorname{dist}(z,\overline{\Omega}_c)
\right\}
\leq
\max\{M\varepsilon_c,0\}
=M\varepsilon_c
\end{align}
for every $c\geq c_0$, where $M$
does not depend on $c$--which establishes (\ref{E:Lipschitz-estimate}).

\medskip

Now, we choose a point $q\in\widehat{K}\setminus K.$ Since $K$ is compact and $q\notin K$,
we have
\begin{align}\label{E: positive-distancefrom K}
d:=\operatorname{dist}(q,K)>0.
\end{align}
\noindent We shall show that $q\notin\overline{\Omega}_c$ for every sufficiently large $c$. Let $w\in\overline{\Omega}_c$.
By (\ref{eq:step4-one-sided}),
there exists $z_w\in K$ such that $\|w-z_w\|\leq M\varepsilon_c.$ By the triangle inequality, $\|q-w\|
\geq
\|q-z_w\|-\|z_w-w\|.$ Since $z_w\in K$, by (\ref{E: positive-distancefrom K}),
$\|q-z_w\|\geq d.$ Consequently, $\|q-w\|
\geq d-M\varepsilon_c.$ Taking the infimum over all
$w\in\overline{\Omega}_c$ gives
\begin{align}\label{eq:step4-exterior-distance}
\operatorname{dist}(q,\overline{\Omega}_c)
\geq d-M\varepsilon_c.
\end{align}

Because $\varepsilon_c\to0$,
there exists $c_1\geq c_0$ such that $M\varepsilon_c<\frac d2$ for $c\geq c_1$. For every such $c$,
\eqref{eq:step4-exterior-distance}
implies $\operatorname{dist}(q,\overline{\Omega}_c)
\geq\frac d2>0.$ In particular, $q\notin\overline{\Omega}_c$ for $
c\geq c_1.$ Combining these facts, we obtain $q\in\widehat{K}\setminus\overline{\Omega}.$ 
This completes the proof.

\medskip

We next show that $\Omega_c$ is not Runge in $\mathbb{C}^2$
for all sufficiently large $c$. Motivated by Wermer's
classical construction \cite{Wermer1980}, we move the
boundary circle of $K$ slightly into $\Omega_c$ and
show that its polynomial hull contains a point outside
$\overline{\Omega}_c$.

\begin{proposition}\label{prop:not-runge}
For all sufficiently large $c$, the domain
$\Omega_c=\Phi(D_c)$ is not Runge in $\mathbb{C}^2$.
\end{proposition}

\begin{proof}
Fix $c>0$ sufficiently large so that $q:=\left(\frac12,0\right)
\notin\overline{\Omega}_c.$ Recall that
\begin{align*}
\rho_c(u,v)
=
y_2+c(y_1^2+y_2^2)
+\chi(x_1^2+x_2^2-1),
\end{align*}
where $u=x_1+iy_1$ and $v=x_2+iy_2$.

\smallskip

For $t\geq0$, consider the curve
\begin{align*}
C_t=
\left\{
\left(a-\frac t2,b-\frac{it}{2}\right):
a^2+b^2=1,\ a,b\in\R
\right\}.
\end{align*}

We first show that $C_t\subset D_c$
for all sufficiently small $t>0$. For $a^2+b^2=1$, define
\begin{align*}
H(a,b,t)
:=
\rho_c\left(a-\frac t2,
b-\frac{it}{2}\right).
\end{align*}

A direct computation gives
\begin{align*}
H(a,b,t)
=
-\frac t2+\frac{ct^2}{4}
+\chi\left(-at+\frac{t^2}{4}\right).
\end{align*}

Since $\chi(0)=\chi'(0)=0$, we have
\begin{align*}
H(a,b,0)=0,
\qquad
\frac{\partial H}{\partial t}(a,b,0)
=-\frac12.
\end{align*}

\noindent By continuity and compactness of the unit
circle, there exists $t_0>0$ such that $\frac{\partial H}{\partial t}(a,b,t)
<-\frac{1}{4}$ for every $a^2+b^2=1$ and $t\in (0,t_0)$. Consequently, $H(a,b,t)<-\frac{t}{4}<0$ for every 
$t \in (0,t_0).$ Thus
$C_t\subset D_c$
for $t\in (0,t_0).$ Now put $\Gamma_t:=\Phi(C_t).$ Since $\Omega_c=\Phi(D_c)$, we have $\Gamma_t\Subset\Omega_c.$ Recall that
\begin{align*}
\Phi(u,v)
=
\bigl(u+iv,F(u+iv,u-iv)\bigr),
\end{align*}
where $F(z,\zeta)
=-(1+i)\zeta+iz\zeta^2+z^2\zeta^3.$
For $(u,v)\in C_t$, we have $u=a-\frac t2,$ and  $v=b-\frac{it}{2},$ where $a,b\in\R$ and $a^2+b^2=1$.
Putting $z=a+ib$, we obtain
\begin{align*}
u+iv&=a+ib=z,\\
u-iv&=a-ib-t=\overline z-t.
\end{align*}
Since $|z|=1$, we have $\overline z=1/z$.
Thus $u+iv=z,$ and $u-iv=\frac1z-t.$ Therefore
\begin{align*}
\Gamma_t=
\left\{
\left(z,F\left(z,\frac1z-t\right)\right):
|z|=1
\right\}.
\end{align*}

By a simple computation, we have
\begin{align*}
F\left(z,\frac1z-t\right)
=-(2+i)t+(3+i)t^2z-t^3z^2.
\end{align*}

Let us define
\begin{align*}
g_t(z): =-(2+i)t+(3+i)t^2z-t^3z^2.
\end{align*}
Then $\Gamma_t=\{(z,g_t(z)):|z|=1\}.$ Since $g_t$ is a holomorphic polynomial,
the maximum modulus principle implies that
the analytic disc $A_t=\{(z,g_t(z)):|z|\leq1\}$ is contained in $\widehat{\Gamma_t}$. In particular,
\begin{align*}
q_t:=
\left(\frac12,g_t\left(\frac12\right)\right)
\in\widehat{\Gamma_t}.
\end{align*}

\noindent Since $g_t\left(\frac12\right)\longrightarrow0$
as $t\to0^+,$ we have $q_t\longrightarrow
\left(\frac12,0\right)=q$ as  $t\to0^+.$ By our choice of $c$,
$q\notin\overline{\Omega_c}$.
Hence $q_t\notin\overline{\Omega_c}$ for all sufficiently small $t>0$.

\smallskip

We have therefore found a compact set
$\Gamma_t\subset\Omega_c$ such that $\widehat{\Gamma_t}
\not\subset\Omega_c.$ Consequently, $\Omega_c$ is not Runge
in $\mathbb{C}^2$. This completes the proof.
\end{proof}

\section{Proof of Theorem~\ref{thm:main-truncated}}

We will prove that $K_R$ is polynomially convex for every $R\in (0,1)$. The proof for $E_R$ is analogous and is omitted.

\medskip

Consider $P(z,w)=zw$. On $K_R$,
\begin{align}\label{E:Image_K_R}
 P(z,f(z))=|z|^2(|z|^2-1)(|z|^2+1+i)   
\end{align}

\noindent Consider the complex-valued function
$S:\mathbb{R}\to\mathbb{C}$ defined by
\begin{align*}
  S(t)=t(t-1)(t+1+i)=t^3+it^2-(1+i)t.  \end{align*}
\noindent {\bf Claim.} For $L\in [0,1),$ $S$ is injective on $[0,L]$.

Let $x,y\in[0,L]$ and suppose $S(x)=S(y)$.
For real $t$,
\begin{align*}
 \operatorname{Re}S(t)=t(t-1)(t+1),\qquad
 \operatorname{Im}S(t)=t(t-1).
\end{align*}
Comparing imaginary parts gives $(x-y)(x+y-1)=0.$ Thus either $x=y$ or $x+y=1$.
Assume $x+y=1$ and $x\ne y$.
If $x,y\in(0,1)$, the common number
$A=x(x-1)=y(y-1)$ is nonzero. Comparing real parts then gives $A(x+1)=A(y+1),$ so $x=y$, a contradiction. If one of $x,y$ equals $0$, the other would have to be $1$,
which is excluded by $L<1$.
This proves injectivity on $[0,L]$.

\smallskip

\noindent This establishes the above claim.\hfill $\blacktriangleleft$

\smallskip

In view of the above claim and \eqref{E:Image_K_R}, we obtain
 $P(K_R)=S([0,R^2])$ is a Jordan arc. Therefore, $\mathcal{P}( P(K_R))=\mathcal{C}(P(K_R))$
\medskip

\noindent Fix $c\in P(K_R)$.
If $c=0$, the corresponding fiber of $K_R$ is
$\{(0,0)\}$. If $c\ne0$, there is a unique
$t\in(0,R^2]$ with $S(t)=c$, and
\begin{align*}
 K_R\cap P^{-1}(c)
 =\{(z,c/z):|z|=\sqrt t\},
\end{align*}
which is polynomially convex.
The theorem is now follows from Result~\ref{R:PolyCnvx_Stolzenberg}

\section{Proof of Theorem~\ref{T: Disc_Levi-flat}}
\label{sec:levi-flat-counterexample}

\begin{proof}
\noindent
We will use the the map $\Phi,$ the non-polynomially convex totally real disc (Wermer disc) $K:=K_1$, the open set $V:=\Phi(U)$ constructed in  the proof of Theorem~\ref{thm:counterexample}.

\medskip

Consider the standard real hypersurface $H:=
\{(u,v)\in\mathbb{C}^2:
\operatorname{Im}v=0\}.$ This is a nonsingular real-analytic
Levi-flat hypersurface in $\mathbb{C}^2$. Now define
\begin{align*}
  M:=\Phi(H\cap U)\subset V.  
\end{align*}

\medskip

\noindent We first show that $M$ is a nonsingular
real-analytic Levi-flat hypersurface
in $V$.

Let $\Psi=(\Psi_1,\Psi_2)
=
(\Phi|_U)^{-1}:V\longrightarrow U$ denote the holomorphic inverse of $\Phi$. Define $\rho:V\longrightarrow\mathbb{R}$ by $\rho(z,w)
=
\operatorname{Im}\Psi_2(z,w).$ By the definition of $M$,
\begin{align*}
M=\{(z,w)\in V:\rho(z,w)=0\}.
\end{align*}

\noindent Since $\Psi_2$ is holomorphic,
$\rho$ is real-analytic and
pluriharmonic on $V$. Therefore, $M$ is real-analytic Lev-flat hypersurface in $V.$ On the other hand, since $\Psi$ is biholomorphic,
its complex differential $D\Psi$
is invertible at every point of $V$. Consequently, $M=\{\rho=0\}$ is a nonsingular.

\medskip

\noindent Next we show that $V$ cannot be a Runge domain
in $\mathbb{C}^2$. Indeed, define $h(w)=\exp(i\Psi_2(w)).$ Then $h$ is holomorphic on $V$ and $M=\{w\in V:|h(w)|=1\}.$ Moreover, $dh(w)=ih(w)d\Psi_2(w)\neq0$ for every $w\in V$, since $\Psi$ is
biholomorphic. If $V$ is Runge, Result~\ref{R:Levi-flat-polynomial-convexity} implies 
 $K$ is
polynomially convex, which is a contradiction. Therefore, $V$ is not Runge in $\mathbb{C}^2.$ 

\medskip
\noindent By (\ref{E:Real Disc}), $\Delta_{\mathbb{R}}\subset H\cap U,$ $K=\Phi(\Delta_{\mathbb{R}})\subset M.$ 

 \medskip
 
Thus, Claim 1 and Claim 2 together complete the proof.

\end{proof}

\section{Proof of Theorem~\ref{thm:main-inward}}

\begin{proof}[Proof of Theorem~\ref{thm:main-inward}]
Let $\{\phi_t\}_{|t|<\varepsilon}$ be the
local holomorphic flow given in the
hypothesis.

\smallskip

We first show that $\widehat{K}\subset\overline\Omega.$ To show this, let $q\in\widehat{K}$ be arbitrary. For each $0<t<\varepsilon$, the map
$\phi_t$ is holomorphic on a neighborhood
of $\widehat{K}$. Hence, by
Lemma~\ref{lem:holomorphic-image-hull}, $\phi_t(q)\in\widehat{\phi_t(K)}.$ By \eqref{eq:inward-flow-condition},
$\phi_t(K)$ is a compact subset of
$\Omega$. Since $\Omega$ is Runge in
$\mathbb{C}^2$, we have $\widehat{\phi_t(K)}\subset\Omega.$ Consequently, $\phi_t(q)\in\Omega
\qquad(0<t<\varepsilon).$
 Since $\phi_0=\operatorname{id}$ and
the flow depends continuously on $t$,
we have $\phi_t(q)\longrightarrow q
\qquad\text{as }t\to0^+.$ Therefore $q\in\overline\Omega$.
As $q\in\widehat{K}$ was arbitrary,
we conclude that $\widehat{K}\subset\overline\Omega.$

\smallskip

The result is now follows from Lemma~\ref{lem:hull-containment}.
\end{proof}

\section{Proof of Theorem~\ref{T:vector-field-equivalence}}
First we will prove the following proposition.

\begin{proposition}\label{prop:converse-vector-field}
Let $K\subset\mathbb{C}^2$ be a polynomially convex,
$C^2$-smooth embedded totally real closed disc,
and let $\rho$ be a real-valued $C^1$-smooth
function defined on a neighborhood of $K$,
satisfying $d\rho_p\neq0
\qquad\text{for every }p\in K.$ Then there exists a polynomial holomorphic
vector field $X:\mathbb{C}^2\to\mathbb{C}^2$ such that $d\rho_p(X(p))<0
\qquad\text{for every }p\in K.$ Moreover, if $\rho|_K=0$, then the local
holomorphic flow $\{\phi_t\}$ generated by $X$
satisfies $\phi_t(K)\subset\{\rho<0\}$ for all sufficiently small $t>0$.
\end{proposition}

\begin{proof}
Since $K$ is polynomially convex and totally real,
Result~\ref{OPW1983} gives
$\mathcal P(K)=\mathcal C(K)$. Define a continuous ambient vector field
$N:K\to\C^2$ by $N(p)=-
\bigl(\rho_{\bar z_1}(p),
      \rho_{\bar z_2}(p)\bigr).$ Under the identification $\C^2\cong\R^4$,
this vector field corresponds to
$-\frac12\nabla_{\R}\rho$.
Since $\rho$ is real-valued, we have
\begin{align*}
d\rho_p(N(p))
=-2\sum_{j=1}^2
|\rho_{z_j}(p)|^2
=-\frac12\|\nabla_{\R}\rho(p)\|^2<0
\end{align*}
for every $p\in K$, where the last inequality
follows from $d\rho_p\neq0$. By compactness of $K$ there exists $\eta>0$
such that $d\rho_p(N(p))
\leq-\eta
$ for all $p\in K.$ Set $L:=\sup_{p\in K}\|d\rho_p\|.$ Since $\mathcal P(K)=\mathcal C(K)$,
we can approximate the components of $N$
uniformly on $K$ by holomorphic polynomials
$P_1,P_2$ such that the polynomial vector
field $X=(P_1,P_2)$ satisfies
\begin{align*}
\sup_{p\in K}\|X(p)-N(p)\|
<\frac{\eta}{2L}.
\end{align*}
Consequently, for every $p\in K$,
\begin{align*}
d\rho_p(X(p))
&=d\rho_p(N(p))
  +d\rho_p(X(p)-N(p))
\leq-\eta+L\|X(p)-N(p)\|
<-\frac{\eta}{2}<0.
\end{align*}

Now suppose that $\rho|_K=0$, and let $X$
be the polynomial holomorphic vector field
constructed above. Since $X$ is entire and
$K$ is compact, the local existence theorem
for holomorphic ordinary differential equations
(see Result~\ref{res:local-holomorphic-flow})
gives an open neighborhood $U$ of $K$
and $\varepsilon>0$ such that the local
holomorphic flow $\{\phi_t\}_{|t|<\varepsilon}$
generated by $X$ is defined on $U$ and satisfies
\begin{align*}
 \frac{d}{dt}\phi_t(z)=X(\phi_t(z)),
\qquad
\phi_0(z)=z.   
\end{align*}

\smallskip

\noindent For each $p\in K$, set $h_p(t):=\rho(\phi_t(p)).$ By the chain rule, $h_p'(t)
=d\rho_{\phi_t(p)}
\bigl(X(\phi_t(p))\bigr).$ In particular, $h_p'(0)=d\rho_p(X(p))
\leq-\frac{\eta}{2}.$ By continuity and compactness of $K$,
there exists $t_0>0$ such that $h_p'(t)\leq-\frac{\eta}{4}
\qquad(p\in K,\ 0\leq t<t_0).$ Since $h_p(0)=\rho(p)=0$, we obtain
\begin{align*}
\rho(\phi_t(p))
=\int_0^t h_p'(s)\,ds
\leq-\frac{\eta t}{4}<0
\end{align*}
for every $p\in K$ and $0<t<t_0$.
Therefore, $\phi_t(K)\subset\{\rho<0\}
\qquad(0<t<t_0).$
\end{proof}

\begin{proof}[Proof of Theorem~ \ref{T:vector-field-equivalence}]

\noindent {\bf {$(i) \implies (iii)$}}: This follows from Proposition~\ref{prop:converse-vector-field}.

\smallskip

\noindent {\bf {$(iii) \implies (ii)$}}: This implication is immediate.

\smallskip

\noindent {\bf {$(ii) \implies (iv)$}}: Let $X$ be an entire holomorphic vector
field satisfying $d\rho_p(X(p))<0$ for
$p\in K,$ and let $\{\phi_t\}$ denote its local
holomorphic flow. Since $X$ is entire and $\widehat{K}$
is compact, by Result~\ref{res:local-holomorphic-flow}, there exist a neighborhood
$U$ of $\widehat{K}$ and $\varepsilon>0$
such that $\phi_t$ is defined and
holomorphic on $U$ for every
$|t|<\varepsilon$. By compactness of $K$, there exists
$\eta>0$ such that $d\rho_p(X(p))\leq-\eta$ for every $p\in K.$ For each $p\in K$, define $h_p(t):=\rho(\phi_t(p)).$ Then $h_p'(0)=d\rho_p(X(p))\leq-\eta.$ By continuity and compactness, after
decreasing $\varepsilon$ if necessary, $h_p'(t)\leq-\frac{\eta}{2}
$ for $p\in K,$ and $ t\in (0,\varepsilon).$
Consequently,
\begin{align*}
\rho(\phi_t(p))
=\int_0^t h_p'(s)\,ds
\leq-\frac{\eta t}{2}<0.
\end{align*}
Hence $\phi_t(K)\subset\Omega$ for 
for $t\in (0,\varepsilon),$ proving (iv).

\smallskip

\noindent\textbf{$(iv)\implies(i)$}:
This follows from Theorem~\ref{thm:main-inward}.

\smallskip

\noindent\textbf{$(i)\implies(v)$}:
This follows from Proposition~\ref{prop:converse-vector-field}.

\smallskip

\noindent\textbf{$(v)\implies(iv)$}:
This implication is immediate.

\smallskip

\noindent
Thus, all five statements are equivalent.
\end{proof}

Theorem~\ref{thm:counterexample} answers
Question~\ref{Q: General case} negatively.
However, the constructed domain is not Runge in $\mathbb{C}^2$.
This raises the question of whether the conclusion
holds for Runge domains.

\begin{question}\label{Q:general-runge}
Let $\Omega\Subset\mathbb{C}^2$ be a bounded Runge domain with
$C^\infty$-smooth strongly pseudoconvex boundary, and let
$K\subset\partial\Omega$ be a smoothly embedded totally real
closed disc. Must $\widehat{K}\subset\overline\Omega$ hold?
\end{question}

\begin{remark}
If $\widehat{K}\subset\overline\Omega$, then the result
of Forstneri\v{c} and Stout \cite{ForSto1991}, together
with the Oka--Weil approximation theorem, implies that
$\widehat{K}=K$. Thus, an affirmative answer to the above question would guarantee the
polynomial convexity of $K$.
\end{remark}

\noindent {\bf Acknowledgements.}  The author acknowledges the financial support of the NBHM, DAE, Government of India, under the NBHM Postdoctoral Fellowship (Reference No. 0204/16(4)/2024/R$\&$D-II/6748, dated 09.05.2024).

\bibliographystyle{plain}
\bibliography{biblio.bib}

\end{document}